\documentclass[proceedings%
]{aofa}

\usepackage[utf8]{inputenc}
\usepackage{subfigure}

\usepackage[round]{natbib}

\author[S. Franceschi, I. Kurkova and K. Raschel]{Sandro Franceschi\addressmark{1,2}\and Irina Kurkova\addressmark{1}\and
  Kilian Raschel\addressmark{2}}

\title[Analytic approach for reflected Brownian motion in the quadrant]{Analytic approach for reflected Brownian motion in the quadrant}

\address{\addressmark{1}Laboratoire de Probabilit\'es et
        Mod\`eles Al\'eatoires, Universit\'e Pierre et Marie Curie,
        4 Place Jussieu, 75252 Paris Cedex 05, France\\
\addressmark{2}Laboratoire de Math\'ematiques et Physique Th\'eorique, Universit\'e de Tours, Parc de Grandmont, 37200 Tours, France}
\keywords{Reflected Brownian motion in the quadrant; (Random) Walks in the quarter plane; Stationary distribution; Laplace transform; Generating function; Boundary value problem; Asymptotic analysis}

\usepackage{epic,pstricks,amsmath,enumitem,tikz}
\usepackage{amsfonts}
\usepackage{amssymb}

\usepackage[mathscr]{euscript}
\colorlet{darkgreen}{green!50!black}

\newtheorem{thm}{Theorem}

\newtheorem{lem}[thm]{Lemma}

\renewcommand{\geq}{\geqslant}
\renewcommand{\leq}{\leqslant}

\usepackage{dsfont}

\usepackage{enumitem}

\usepackage{eqnarray,cases}

\usepackage{float}

\newcommand{\s}{\mathscr{S}}

\newcommand{\G}{\mathscr G_\mathscr R}
\newcommand{\R}{\mathscr{R}}

\begin{document}
\maketitle
\begin{abstract}
\paragraph{Abstract.}
Random walks in the quarter plane are an important object both of combinatorics and probability theory. Of particular interest for their study, there is an analytic approach initiated by Fayolle, Iasnogorodski and Maly\v sev, and further developed by the last two authors of this note. 
The outcomes of this method are explicit expressions for the generating functions of interest, asymptotic analysis of their coefficients, etc. Although there is an important literature on reflected Brownian motion in the quarter plane (the continuous counterpart of quadrant random walks), an analogue of the analytic approach has not been fully developed to that context. The aim of this note is twofold: it is first an extended abstract of two recent articles of the authors of this paper
, which propose such an approach; we further compare various aspects of the discrete and continuous analytic approaches.
 \end{abstract}

\section{Introduction}
\label{intro}

\subsection{Random walks in the quarter plane}
Since the seventies and the pioneered papers \cite{Ma72,FaIa79}, random walks in the quarter plane (cf.\ Figure \ref{fig:marches_avec_poids}) are extensively studied. They are indeed an important object of probability theory and have been studied for their recurrence/transience, for their links with queueing systems (\cite{FaIa79}), representation theory (\cite{Bi92}), potential theory. Moreover, the state space $\mathbb N^2=\{0,1,2,\ldots \}^2$ offers a natural framework for studying any two-dimensional population; accordingly, quadrant walks appear as models in biology and in finance (\cite{CoLa13}). Another interest of random walks in the quarter plane is that in the large class of random processes in cones, they form a family for which remarkable exact formulas exist. Moreover, quadrant walks are popular in combinatorics, see \cite{BMMi10,BoKa10,KuRa12}. Indeed, many models of walks are in bijection with other combinatorial objects: maps, permutations, trees, Young tableaux, etc. In combinatorics again, famous models have emerged from quadrant walks, as Kreweras' or Gessel's ones, see \cite{BMMi10,BoKa10}. Finally, walks in the quarter plane are interesting for the numerous tools used for their analysis: combinatorial (\cite{BMMi10}), from complex analysis (\cite{Ma72,FaIa79,FaIaMa99,KuRa11,KuRa12,BeBMRa15}), computer algebra (\cite{BoKa10}), for instance.

\subsection{Issues and technicalities of the analytic approach}

In the literature (see, e.g., \cite{Ma72,FaIa79,FaIaMa99,KuRa11,KuRa12}), the analytic approach relies on six key steps: 
\begin{enumerate}[label={\rm (\roman{*})},ref={\rm (\roman{*})}]
     \item\label{enumi:functional_equation}Finding a functional equation between the generating functions of interest;
     \item\label{enumi:BVP_1}Reducing the functional equation to boundary value problems (BVP);
     \item\label{enumi:BVP_2}Solving the BVP;
     \item\label{enumi:group}Introducing the group of the walk;
     \item\label{enumi:Riemann}Defining the Riemann surface naturally associated with the model, continuing meromorphically the generating functions and finding the conformal gluing function;
     \item\label{enumi:asymptotic_technique}Deriving the asymptotics of the (multivariate) coefficients.
\end{enumerate}
Before commenting these different steps, let us note that altogether, they allow for studying the following three main problems:
\begin{enumerate}[label={\rm (P\arabic{*})},ref={\rm (P\arabic{*})}]
     \item\label{enumi:explicit}Explicit expression for the generating functions of interest (needs \ref{enumi:functional_equation}, \ref{enumi:BVP_1}, \ref{enumi:BVP_2} and \ref{enumi:Riemann});
     \item\label{enumi:algebraic}Algebraic nature of these functions (needs \ref{enumi:group} and \ref{enumi:Riemann});
     \item\label{enumi:asymptotic}Asymptotics of their coefficients in various regimes (needs \ref{enumi:BVP_2}, \ref{enumi:Riemann} and \ref{enumi:asymptotic_technique}).     
\end{enumerate}
The point \ref{enumi:functional_equation} reflects the inherent properties of the model and is easily obtained. Point \ref{enumi:BVP_1}, first shown in \cite{FaIa79}, is now standard (see \cite{FaIaMa99}) and follows from algebraic manipulations of the functional equation of \ref{enumi:functional_equation}. Item \ref{enumi:BVP_2} uses specific literature devoted to BVP (our main reference for BVP is the book of \cite{Li00}). The idea of introducing the group of the model \ref{enumi:group} was proposed in \cite{Ma72}, and brought up to light in the combinatorial context in \cite{BMMi10}. Point \ref{enumi:Riemann} is the most technical a priori; it is however absolutely crucial, as it allows to access key quantities (as a certain conformal gluing function which appears in the exact formulation of \ref{enumi:BVP_2}). Finally, \ref{enumi:asymptotic_technique} uses a double refinement of the classical saddle point method: the uniform steepest descent method.

\subsection{Reflected Brownian motion in the quarter plane}
There is a large literature on reflected Brownian motion in quadrants (and also in orthants, generalization to higher dimension of the quadrant), to be rigorously introduced in Section \ref{sec:BM}. First, it serves as an approximation of large queuing networks (see \cite{Fo84,BaFa87}); this was the initial motivation for its study. In the same vein, it is the continuous counterpart of (random) walks in the quarter plane. In other directions, it is studied for its Lyapunov functions in \cite{DuWi94}, cone points of Brownian motion in \cite{LG87}, intertwining relations and crossing probabilities in \cite{Du04}, and of particular interest for us, for its recurrence/transience in \cite{HoRo93}. The asymptotic behavior of the stationary distribution (when it exists) is now well known, see \cite{HaHa09,DaMi11,FrKu16}. There exist, however, very few results giving exact expressions for the stationary distribution. Let us mention \cite{Fo84} (dealing with the particular case of a Brownian motion with identity covariance matrix), \cite{BaFa87} (on a diffusion having a quite special behavior on the boundary), 
\cite{harrison_multidimensional_1987,DiMo-09} (on the special case when stationary densities are exponential) and \cite{FrRa16} (on the particular case of an orthogonal reflection). We also refer to \cite{BuChMaRa15} for the analysis of reflected Brownian motion in bounded planar domains by complex analysis techniques.

\subsection{Main results and plan}
This note is an extended abstract of the papers \cite{FrKu16,FrRa16}, whose main contributions are precisely to export the analytic method for reflected Brownian motion in the quarter plane. Our study constitutes one of the first attempts to apply these techniques to the continuous setting, after \cite{Fo84,BaFa87}. In addition of reporting about the works \cite{FrKu16,FrRa16}, we also propose a comparative study of the discrete/continuous cases. 

Our paper is organized as follows: Section \ref{sec:RW} concerns random walks and Section \ref{sec:BM} Brownian motion. For clarity of exposition we have given the same structure to Sections \ref{sec:RW} and \ref{sec:BM}: in Section~\ref{subsec:functional_equation_RW}/\ref{subsec:functional_equation_BM} we first state the key functional equation (a kernel equation), which is the starting point of our entire analysis. We study the kernel (a second degree polynomial in two variables). In Section \ref{subsec:BVP_RW}/\ref{subsec:BVP_BM} we state and solve the BVP satisfied by the generating functions. We then move to asymptotic results (Section \ref{subsec:asymp_RW}/\ref{subsec:asymp_BM}). In Section \ref{subsec:Riemann_surface_RW}/\ref{subsec:Riemann_surface_BM} we introduce the Riemann surface of the model and some important related facts. 

\section{Random walks in the quarter plane}
\label{sec:RW}

This section is devoted to the discrete case and is based mainly on \cite{FaIaMa99}.

\subsection{Functional equation}
\label{subsec:functional_equation_RW}

One considers a piecewise homogeneous random walk with sample paths in $\mathbb N^2$. There are four domains of spatial homogeneity (the interior of $\mathbb N^2$, the horizontal and vertical axes, the origin), inside of which the transition probabilities (of unit size) are denoted by $p_{i,j}$, $p_{i,j}'$, $p_{i,j}''$ and $p_{i,j}^0$, respectively. See Figure \ref{fig:marches_avec_poids}. The inventory polynomial of the inner domain is called the kernel and equals
\begin{equation}
\label{eq:kernel_RW}
     K(x,y) =  xy \{ \textstyle\sum_{-1\leq i,j\leq1} p_{i,j} x^i y^j-1 \}.
\end{equation}
The inventory polynomials associated to the other homogeneity domains are
\begin{equation*}  
     k(x,y) =  \textstyle x \{\sum p'_{i,j} x^i y^j - 1\}, \quad \widetilde{k}(x,y) =  y \{\sum p''_{i;j} x^i y^j - 1 \},\quad k_0(x,y) =   \{\sum p^0_{i,j} x^i y^j - 1\}. 
\end{equation*}
Assuming the random walk ergodic (we refer to \cite[Theorem~1.2.1]{FaIaMa99} for necessary and sufficient conditions), we denote the invariant measure by $\{\pi_{i,j}\}_{i,j \geq 0}$ and introduce the generating functions
\begin{equation*}  
     \textstyle\pi(x,y) =  \sum_{i,j \geq 1} \pi_{i,j} x^{i-1} y^{j-1}, \quad\pi(x) =  \sum_{i\geq 1} \pi_{i,0} x^{i-1}, \quad \widetilde{\pi}(y) = \sum_{j\geq 1} \pi_{0,j} y^{j-1}.
\end{equation*}
Writing the balance equations at the generating function level, we have (see \cite[Equation~(1.3.6)]{FaIaMa99} for the original statement):
\begin{lem}
The fundamental functional equation holds 
\begin{equation} 
\label{eq:functional_equation_RW}
      -K(x,y) \pi(x,y) = k(x,y) \pi(x) + \widetilde{k}(x,y)\widetilde{\pi}(y) + k_0(x,y)\pi_{0,0}.
\end{equation}
\end{lem}
Equation \eqref{eq:functional_equation_RW} holds a priori in the region $\{(x,y)\in \mathbb{C}^{2} : \vert x\vert\leq 1,\vert y\vert\leq1\}$. Indeed, the $\pi_{i,j}$ sum up to $1$, so that the generating functions $\pi(x,y)$, $\pi(x)$ and $\widetilde{\pi}(y)$ are well defined on the (bi)disc. The identity \eqref{eq:functional_equation_RW} is a kernel equation, and a crucial role will be played by the kernel \eqref{eq:kernel_RW}. This polynomial $K$ is of second order in both $x$ and $y$; its roots $X(y)$ and $Y(x)$ defined by 
\begin{equation}
\label{eq:definition_branches_RW}
     K(X(y),y)=K(x,Y(x))=0
\end{equation}     
are thus algebraic of degree $2$. Writing the kernel as $K(x,y)= a(y)x^2+ b(y)x+ c(y)$ and defining its discriminant $ d(y)= b(y)^2-4 a(y) c(y)$, one has obviously
\begin{equation*}
     X(y)=\frac{-b(y)\pm\sqrt{ d(y)}}{2 a(y)}.
\end{equation*}
The polynomial $d$ has three or four roots, and exactly two of them are located in the unit disc, see \cite[Lemma 2.3.8]{FaIaMa99}. They are named $y_1,y_2$, cf.\ Figure \ref{fig:BVP_RW}. On $(y_1,y_2)$ one has $d(y)<0$, so that the two values (or branches) of $X(y)$ (that we shall call $X_0(y)$ and $X_1(y)$) are complex conjugate of one another. In particular, the set
\begin{equation*}
     \mathscr M=X([y_1,y_2])=\{x\in\mathbb C: K(x,y)=0\text{ and } y\in[y_1,y_2]\}
\end{equation*}
is symmetrical w.r.t.\ the real axis (Figure \ref{fig:BVP_RW}). This curve will be used to set a boundary condition for the unknown function $\pi$ (Lemma \ref{lem:BVP_RW}).

\unitlength=1cm
\begin{figure}[t]        
\begin{center}
        \begin{picture}(6.5,5.5)
    \thicklines
    \put(1,1){{\vector(1,0){4.5}}}
    \put(1,1){\vector(0,1){4.5}}
        \linethickness{0.3mm}
    \put(1,2){\dottedline{0.15}(0,0)(4.5,0)}
    \put(1,3){\dottedline{0.15}(0,0)(4.5,0)}
    \put(1,4){\dottedline{0.15}(0,0)(4.5,0)}
    \put(1,5){\dottedline{0.15}(0,0)(4.5,0)}
    \put(2,1){\dottedline{0.15}(0,0)(0,4.5)}
    \put(3,1){\dottedline{0.15}(0,0)(0,4.5)}
    \put(4,1){\dottedline{0.15}(0,0)(0,4.5)}
    \put(5,1){\dottedline{0.15}(0,0)(0,4.5)}
    \put(4,4){\textcolor{blue}{\vector(1,1){1}}}
    \put(4,4){\textcolor{blue}{\vector(-1,-1){1}}}
    \put(4,4){\textcolor{blue}{\vector(1,0){1}}}
    \put(4,4){\textcolor{blue}{\vector(-1,0){1}}}
    \put(4,4){\textcolor{blue}{\vector(0,1){1}}}
    \put(4,4){\textcolor{blue}{\vector(0,-1){1}}}
    \put(4,4){\textcolor{blue}{\vector(-1,1){1}}}
    \put(4,4){\textcolor{blue}{\vector(1,-1){1}}}
    \put(4,1){\textcolor{darkgreen}{\vector(1,0){1}}}
    \put(4,1){\textcolor{darkgreen}{\vector(-1,0){1}}}
    \put(4,1){\textcolor{darkgreen}{\vector(0,1){1}}}
    \put(4,1){\textcolor{darkgreen}{\vector(-1,1){1}}}
    \put(4,1){\textcolor{darkgreen}{\vector(1,1){1}}}
    \put(1,4){\textcolor{orange}{\vector(1,0){1}}}
    \put(1,4){\textcolor{orange}{\vector(0,-1){1}}}
    \put(1,4){\textcolor{orange}{\vector(0,1){1}}}
    \put(1,4){\textcolor{orange}{\vector(1,-1){1}}}
    \put(1,4){\textcolor{orange}{\vector(1,1){1}}}
    \put(1,1){\textcolor{red}{\vector(0,1){1}}}
    \put(1,1){\textcolor{red}{\vector(1,1){1}}}
    \put(1,1){\textcolor{red}{\vector(1,0){1}}}
    \put(5.05,5.2){\textcolor{blue}{$p_{i,j}$}}
    \put(2.05,5.2){\textcolor{orange}{$p_{i,j}''$}}
    \put(5.05,2.2){\textcolor{darkgreen}{$p_{i,j}'$}}
    \put(2.05,2.2){\textcolor{red}{$p_{i,j}^{0}$}}
\end{picture}
\end{center}
  \vspace{-10mm}
\caption{Transition probabilities of the reflected random walk in the quarter plane, with four domains of spatial homogeneity}
\label{fig:marches_avec_poids}
\end{figure}
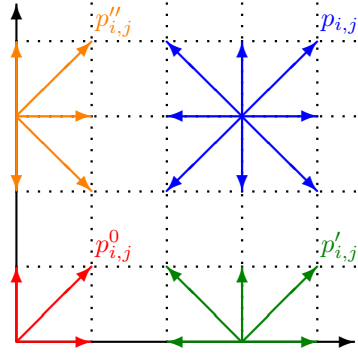

\subsection{Statement and resolution of the BVP}
\label{subsec:BVP_RW}

The analytic approach of \cite{Ma72,FaIa79,FaIaMa99} proposes a way for solving the functional equation \eqref{eq:functional_equation_RW}, by reduction to a BVP. Generally speaking, a BVP consists of a regularity condition and a boundary condition. 


\begin{lem}
\label{lem:BVP_RW}
The function $\pi$ satisfies the following BVP:
\begin{itemize}
     \item $\pi$ is meromorphic in the bounded domain delimitated by $\mathscr M$ and has there identified poles;
     \item for any $x\in\mathscr M$,
     \begin{equation*}
          \frac{k(x,Y_0(x))}{\widetilde k(x,Y_0(x))}\pi(x)-\frac{k(\overline{x},Y_0(\overline{x}))}{\widetilde k(\overline{x},Y_0(\overline{x}))}\pi(\overline{x})=\frac{k_0(\overline{x},Y_0(\overline{x}))}{\widetilde k(\overline{x},Y_0(\overline{x}))}\pi_{0,0}-\frac{k_0(x,Y_0(x))}{\widetilde k(x,Y_0(x))}\pi_{0,0}.
     \end{equation*}     
\end{itemize}     
\end{lem}

\begin{proof}
The regularity condition follows from Theorem \ref{thm:continuation_RW}, which provides a (maximal) meromorphic continuation of the function $\pi$.
We now turn to the boundary condition. For $i\in\{0,1\}$, we evaluate the functional equation \eqref{eq:functional_equation_RW} at $(X_i(y),y)$ and divide by $\widetilde k(X_i(y),y)$. We then make the difference of the identities corresponding to $i=0$ and $i=1$. Finally, we substitute $X_0(y)=x$ and $X_1(y)=\overline{x}$, noting that when $y\in[y_1,y_2]$, $x\in \mathscr M$ by construction. Notice that we have chosen the segment $[y_1,y_2]$ connecting the points inside of the unit disc (Figure \ref{fig:BVP_RW}), in which we know that the generating function $\widetilde \pi$ is well defined. 
\end{proof}

\unitlength=0.6cm
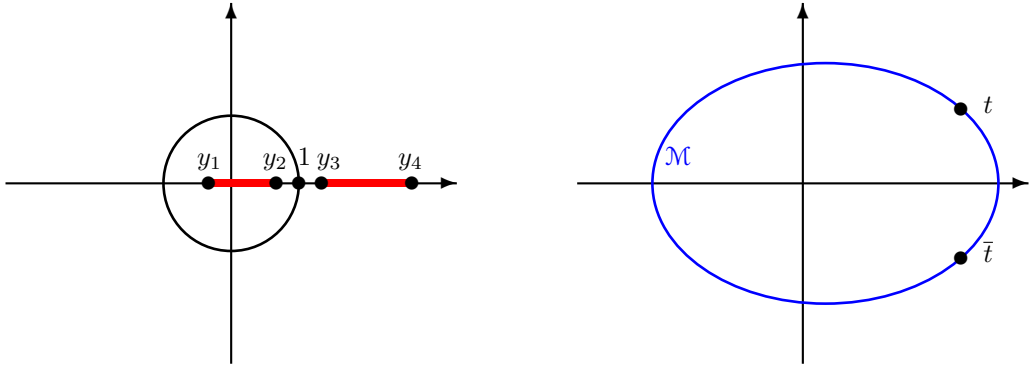
\begin{figure}[h!]
\vspace{15mm}
   \hspace{25mm}
    \begin{tikzpicture}(0,5)
    \thicklines
    \put(-5,0){\vector(1,0){10}}
    \put(0,-4){\vector(0,1){8}}
   \textcolor{black}{\draw[line width=1pt,fill=none] (0,0) ellipse (0.9cm and 0.9cm);}
   \put(1.48,0.4){$1$}
   \put(-0.75,0.4){$y_1$}
   \put(0.70,0.4){$y_2$}
   \put(1.9,0.4){$y_3$}
   \put(3.7,0.4){$y_4$}
    {\put(1.5,0){\textcolor{black}{\circle*{0.3}}}}
    
\linethickness{1mm}
{\put(-0.5,0){\textcolor{red}{\line(1,0){1.5}}}}
{\put(2,0){\textcolor{red}{\line(1,0){2}}}}
{\put(-0.5,0){\textcolor{black}{\circle*{0.3}}}}
{\put(1,0){\textcolor{black}{\circle*{0.3}}}}
{\put(2,0){\textcolor{black}{\circle*{0.3}}}}
{\put(4,0){\textcolor{black}{\circle*{0.3}}}}
     \end{tikzpicture}
     
     \vspace{-25.7mm}
     
   \hspace{90mm}   \begin{tikzpicture}(0,0)
    \thicklines
    \put(-5,0){\vector(1,0){10}}
    \put(0,-4){\vector(0,1){8}}
    \put(-3.05,0.4){\textcolor{blue}{$\mathscr M$}}
   \textcolor{blue}{\draw[line width=1pt,fill=none] (0.3,0) ellipse (2.3cm and 1.6cm);}
   {\put(4.0,1.55){\textcolor{black}{$t$}}}
   {\put(4.0,-1.75){\textcolor{black}{$\overline{t}$}}}
       {\put(3.5,1.65){\textcolor{black}{\circle*{0.3}}}}
       {\put(3.5,-1.65){\textcolor{black}{\circle*{0.3}}}}
     \end{tikzpicture}
\vspace{7mm}
\caption{Left: the polynomial $d$ has three or four roots, denoted by $y_1,y_2,y_3,y_4$; exactly two of them are inside the unit disc. Right: the curve $\mathscr M=X([y_1,y_2])$ is symmetrical w.r.t.\ the horizontal axis}
\label{fig:BVP_RW}
\end{figure}

Although Lemma \ref{lem:BVP_RW} could be written more precisely (by giving the number and the location of the poles of $\pi$), we shall prefer the above version, since we focus in this note on the main ideas of the analytic approach. 

Lemma \ref{lem:BVP_RW} happens to characterize the generating functions, as it eventually leads to an explicit expression for $\pi$, see Theorem \ref{thm:explicit_BVP_RW}. Before stating this central result (borrowed from \cite[Theorem~5.4.3]{FaIaMa99}), we need to introduce a function $w$ called a conformal gluing function. By definition it satisfies $w(x)=w(\overline{x})$ for $x\in\mathscr M$ and is one-to-one inside of $\mathscr M$. This function will be constructed in Theorem \ref{thm:formula_w_RW} of Section \ref{subsec:Riemann_surface_RW}.

\begin{thm}
\label{thm:explicit_BVP_RW}
There exist two functions $f$ and $g$, constructed from $k$, $\widetilde k$, $k_0$ and $w$, such that the following integral formulation for $\pi$ holds:
\begin{equation*}
     \pi(x)=f(x)\int_\mathscr M g(u)\frac{w'(u)}{w(u)-w(x)}\textnormal{d}u.
\end{equation*}
\end{thm}
A similar contour integral representation exists for $\widetilde \pi$, and eventually the functional equation \eqref{eq:functional_equation_RW} provides us with an explicit expression for the bivariate function $\pi(x,y)$.

\subsection{Asymptotics of the stationary probabilities}
\label{subsec:asymp_RW}

  The asymptotics of coefficients $\pi_{i,j}$ of unknown 
  generating functions satisfying the functional equation \eqref{eq:functional_equation_RW} has been obtained by  \cite{Ma73} via analytic arguments.
     He computed the asymptotics of the stationary probabilities 
   $\pi_{i,j}$
       as $i,j \to \infty$ and $j/i = \tan \alpha$, for any given $\alpha \in (0, \pi/2)$. 
   Let us briefly  present these results.
   It is assumed in \cite{Ma73} that the random walk is simple, meaning that 
  \begin{equation}
p_{-1,1}=p_{1,1}=p_{-1,-1}=p_{1,-1}=0.
\label{assimple}
\end{equation}
  It is also assumed that both coordinates of the interior drift vector are negative (as in Figure \ref{fig:drift_reflection}). For $\alpha \in (0,\pi/2)$, we define the point $(x(\alpha), y(\alpha))$ as follows. Introducing  as in \cite{KuRa11} the function $P(u,v)=\sum_{i,j}p_{i,j}e^{iu}e^{jv}$ on ${\mathbb R}^2$, the mapping 
\begin{equation*}
     (u,v) \mapsto \frac{\nabla P(u,v)}{\vert P(u,v)\vert}
\end{equation*}
is a homeomorphism between $\{(u,v) \in {\mathbb R}^2 : P(u,v)=1\}$ and the unit circle. The point $(u(\alpha), v(\alpha))$ is the unique solution to $\frac{\nabla P(u,v)}{\vert P(u,v)\vert}=(\cos \alpha, \sin \alpha)$. Finally, $(x(\alpha),y(\alpha))= (e^{u(\alpha)},e^{v(\alpha)})$.
    
  
     Following \cite{Ma73}, we introduce the sets of parameters 
\begin{eqnarray}
 \mathscr{P}_{--}&=&\big\{(\{p_{i,j}\}, \alpha): k(\psi(x(\alpha), y(\alpha)))\leq 0 \text{ and } \widetilde k(\phi(x(\alpha), y(\alpha)))\leq 0 \big\},\nonumber\\
\mathscr{P}_{+-}&=&\big\{(\{p_{i,j}\},  \alpha): k(\psi(x(\alpha), y(\alpha)))> 0 \text{ and }  \widetilde k(\phi(x(\alpha), y(\alpha)))\leq 0 \big\},\nonumber
\end{eqnarray}
   and $\mathscr{P}_{-+}$ and $\mathscr{P}_{++}$ accordingly.
  The automorphisms $\psi$ and $\phi$ are defined in Section \ref{subsec:Riemann_surface_RW} by \eqref{psiphi}.
   The following theorem is proven in \cite{Ma73}.

\begin{thm}
   Let $(i,j)=(r \cos \alpha,r \sin \alpha)$ with $\alpha \in (0, \pi/2)$. 
Then as $r \to \infty$ we have 
\begin{equation}
\label{eqass}
\pi_{i,j} =(1+o(1))\cdot 
   \left\{
\begin{array}{ll}
\frac{C_0(\alpha) }{\sqrt{r}}x^{-i}(\alpha) y^{-j}(\alpha)\  &\hbox{ in } \mathscr{P}_{--},\\
C_1 p_1^{-i} q_1^{-j}\  &\hbox{ in } \mathscr{P}_{-+},\\
C_2 p_2^{-i} q_2^{-j} & \hbox{ in } \mathscr{P}_{+-},\\
C_1 p_1^{-i} q_1^{-j}  +
C_2 p_2^{-i} q_2^{-j}  &\hbox{ in } \mathscr{P}_{++},
\end{array}
\right.
\end{equation}
   where $C_0$, $C_1$ and $C_2$ are constants that can be expressed in terms of the functions $\pi$ and $\widetilde \pi$. 
The point $(p_1, q_1)$ is a  solution of the system $\{K(x,y)=0,\ k(\psi(x,y))=0\}$ and similarly, $(p_2,q_2)$ is a solution of $\{K(x,y)=0,\ \widetilde k(\phi(x,y))=0\}$.
\end{thm}
   
\begin{proof}  
The stationary probabilities $\pi_{i,j}$ are first written as two-dimensional Cauchy integrals, then reduced via the residue theorem to one-dimensional integrals along some contours. The asymptotics of these integrals is characterized either by the saddle point $(x(\alpha), y(\alpha))$ in the case of the set of parameters $\mathscr{P}_{--}$ or by a pole $(p_1, q_1)$ or $(p_2, q_2)$ that is encountered when moving the integration contour to the saddle point; this happens for the sets of parameters $\mathscr{P}_{-+}$, $\mathscr{P}_{+-}$ and $\mathscr{P}_{++}$. 
\end{proof}

This approach has been applied for the analysis of the join the shortest queue problem in \cite{KuSu03}, and for computing the asymptotics of  Green functions for transient random walks in the quarter plane reflected at the axes (see \cite{KuMa98}) or killed at the axes (cf.\ \cite{KuRa11}). Moreover, as illustrated in \cite{KuRa11,KuSu03}, the assumption \eqref{assimple} is not  crucial for the applicability of the method. The limiting cases  $\alpha=0$ and $\alpha=\pi/2$ can also be treated via this approach, with some additional technical details (the saddle point then coincides with a branch point of the integrand), it is done in \cite{KuRa11}.

\subsection{Riemann surface and related facts}
\label{subsec:Riemann_surface_RW}
In Section \ref{subsec:functional_equation_RW} the set
\begin{equation*}
     \mathscr K=\{(x,y)\in\mathbb C^2: K(x,y)=0\}=\{(x,y)\in\mathbb C^2: \textstyle \sum p_{i,j}x^{i}y^{j}=1\}
\end{equation*}
has appeared very naturally, since in order to state the BVP (our Lemma \ref{lem:BVP_RW}), we introduced the functions $X(y)$ and $Y(x)$, which by construction cancel the kernel, see \eqref{eq:definition_branches_RW}. 

In this section the central idea is to consider the (global) complex structure of $\mathscr K$. The set $\mathscr K$ turns out to be a Riemann surface of genus $1$, i.e., a torus. This simply comes from the reformulation of the identity $K(x,y)=0$ as
\begin{equation*}
     \{2a(y)x+b(y)\}^2=d(y).
\end{equation*}
Moreover, the Riemann surface of the square root of a polynomial of degree $3$ or $4$ is classically a torus (with this terminology, the roots of the discriminant are branch points).

This new point of view on $\mathscr K$ brings powerful tools. Of particular interest is a parametrization of $\mathscr K$ in terms of Weierstrass elliptic functions:
\begin{equation}
\label{eq:parametrization}
     \mathscr K=\{(x(\omega),y(\omega)): \omega\in\mathbb C/(\omega_1\mathbb Z+\omega_2\mathbb Z)\}.
\end{equation}
This parametrization is totally explicit: $x(\omega)$ and $y(\omega)$ are rational functions in the $\wp$-Weierstrass function and its derivative $\wp'$ (see \cite[Lemma 3.3.1]{FaIaMa99}); the periods $\omega_1$ and $\omega_2$ admit expressions as elliptic integrals in terms of $\{p_{i,j}\}$ (cf.\ \cite[Lemma 3.3.2]{FaIaMa99}), etc. Moreover, as any functions of $x$ and/or $y$, the functions $\pi(x)$ and $\widetilde\pi(y)$ can be lifted on $\mathscr K$ by setting 
\begin{equation}
\label{eq:lifted_generating_functions_RW}
     \Pi(\omega)=\pi(x(\omega)),\qquad \widetilde\Pi(\omega)=\widetilde\pi(y(\omega)).
\end{equation}

\subsubsection*{Group of the walk}
Introduced in \cite{Ma72} in a probabilistic context and further used in \cite{FaIaMa99,BMMi10}, the group of the walk is a dihedral group generated by 
\begin{equation}
         \zeta(x,y)=\left(x,\frac{\sum_{i}p_{i,-1}x^{i}}
          {\sum_{i}p_{i,+1}x^{i}}\frac{1}{y}\right),
          \qquad \eta(x,y)=\left(\frac{\sum_{j}p_{-1,j}y^{j}}
          {\sum_{j}p_{+1,j}y^{j}}\frac{1}{x},y\right).
          \label{psiphi}
\end{equation}
(One easily verifies that these generators are idempotent: $\zeta^2=\eta^2=1$.) The group $\langle\zeta,\eta\rangle$ can be finite or infinite, according to the order of the element $\zeta\circ\eta$. The generator $\zeta$ (resp.\ $\eta$) exchanges the roots in $y$ (resp.\ in $x$) of $K(x,y)=0$. Viewed as a group of birational transformations in \cite{BMMi10}, we shall rather see it as a group of automorphisms of the Riemann surface $\mathscr K$. 

This group has many applications. First, it allows for a continuation of the functions $\pi(x)$ and $\widetilde\pi(y)$ (Theorem \ref{thm:continuation_RW} below). It further connects the algebraic nature of the generating functions to the (in)finiteness of the group (Theorem \ref{thm:nature_solution_RW}). Finally, in the finite group case, elementary algebraic manipulations of the functional equations can be performed (typically, via the computation of certain orbit-sums) so as to eventually obtain D-finite expressions for the unknowns, see \cite[Chapter 4]{FaIaMa99} and \cite{BMMi10}.

Using the structure of the automorphisms of a torus, the lifted versions of $\zeta$ and $\eta$ admit simple expressions \cite[Section 3.1.2]{FaIaMa99}: 
\begin{equation}
\label{eq:lifted_expression_automorphisms}
     \zeta(\omega)=-\omega+\omega_1+\omega_2,\qquad \eta(\omega)=-\omega+\omega_1+\omega_2+\omega_3, 
\end{equation}
where, as the periods $\omega_1$ and $\omega_2$, $\omega_3\in(0,\omega_2)$ is an elliptic integral \cite[Lemma 3.3.3]{FaIaMa99}. Accordingly, the group is finite if and only if $\omega_2/\omega_3\in\mathbb Q$, which provides a nice criterion in terms of elliptic integrals.

\subsubsection*{Continuation}
While the generating function $\pi(x)$ is defined through its power series in the unit disc, it is a priori unclear how to continue it to a larger domain. This is however crucial, since the curve $\mathscr M$ on which it satisfies a BVP (Lemma \ref{lem:BVP_RW}) is not included in the unit disc in general.

\begin{thm}
\label{thm:continuation_RW}
The function $\pi$ can be continued as a meromorphic function to $\mathbb C\setminus [x_3,x_4]$.
\end{thm}
We notice that $\mathscr M$ does not intersect $[x_3,x_4]$ by \cite[Theorem 5.3.3]{FaIaMa99}, so that Theorem~\ref{thm:continuation_RW} indeed provides a continuation of the generating function in the domain delimitated by $\mathscr M$.

\begin{proof}
This result, stated as Theorem 3.2.3 in \cite{FaIaMa99}, is a consequence of a continuation of the lifted generating functions \eqref{eq:lifted_generating_functions_RW} on the Riemann surface (or, better, on its universal covering --- but we shall not go into these details here). The continuation on $\mathscr K$ uses the (lifted) functional equation \eqref{eq:functional_equation_RW} and the group of the walk $\langle\zeta,\eta\rangle$.
\end{proof}

\subsubsection*{Conformal mapping}

In the integral expression of Theorem \ref{thm:explicit_BVP_RW}, the conformal gluing function $w$ is all-present, as it appears in the integrand and in $f$ and $g$ as well. The introduction of the Riemann surface $\mathscr K$ allows to derive an expression for this function (this is another major interest of introducing $\mathscr K$). Let us recall that $\omega_1$ and $\omega_2$ are the periods of the elliptic functions of the parametrization \eqref{eq:parametrization}, while $\omega_3$ comes out in the lifted expression \eqref{eq:lifted_expression_automorphisms} of the automorphisms.
\begin{thm}
\label{thm:formula_w_RW}
The conformal gluing function $w$ admits the expression:
\begin{equation*}
     w(x)=\wp(\wp^{-1}(x;\omega_1,\omega_2);\omega_1,\omega_3),
\end{equation*}
where for $i\in\{2,3\}$, $\wp(\, \cdot\, ;\omega_1,\omega_i)$ is the $\wp$-Weierstrass elliptic function with periods $\omega_1$ and $\omega_i$.
\end{thm}
\begin{proof}
While on the complex plane, $\mathscr M$ is a quartic curve, it becomes on $\mathscr K$ much simpler (typically, a segment). This remark (which again illustrates all the benefit of having introduced the Riemann surface) is used in \cite[Section 5.5.2]{FaIaMa99} so as to obtain the above expression for $w$.
\end{proof}

\subsubsection*{Algebraic nature of the generating functions}
Recall that a function of one variable is D-finite if it satisfies a linear differential equation with polynomial coefficients.

\begin{thm}
\label{thm:nature_solution_RW}
If the group is finite, the generating functions $\pi(x)$ and $\widetilde\pi(y)$ are D-finite.
\end{thm}
\begin{proof}
This follows from manipulations on the Riemann surface, see \cite[Chapter~4]{FaIaMa99}. The D-finiteness is proved on $\mathbb R$; refined results (in the combinatorial context of the enumeration of paths) can be found in \cite{BMMi10}, where the D-finiteness is proved on $\mathbb Q$.
\end{proof}

The converse of Theorem \ref{thm:nature_solution_RW} is not shown in full generality. It is true in combinatorics, see \cite{KuRa12}.

\section{Reflected Brownian motion in the quadrant}
\label{sec:BM}
\let\phi=\varphi

\subsubsection*{Defining reflected Brownian motion in the quadrant}

\begin{figure}[hbtp]
\centering
\includegraphics[scale=0.6]{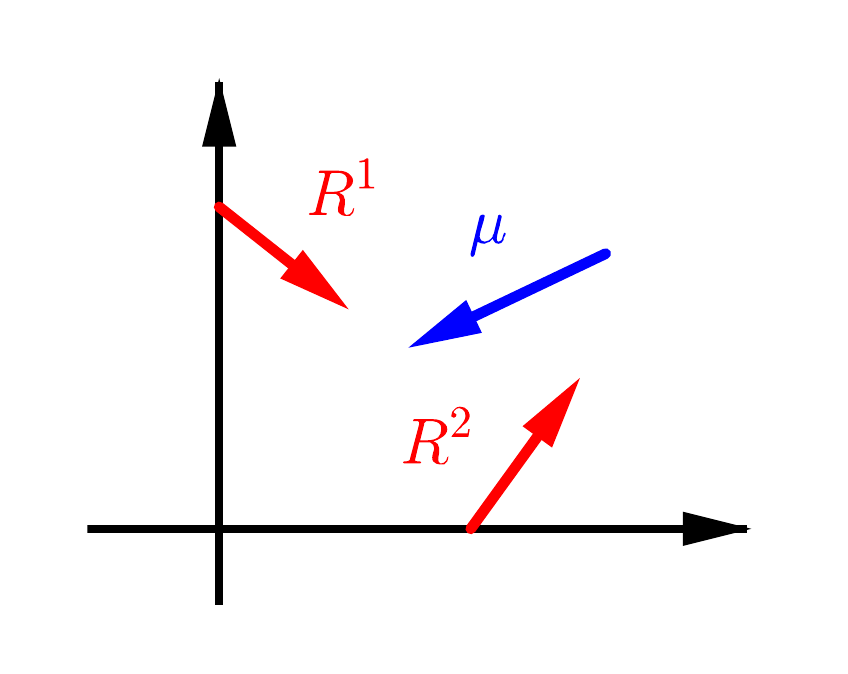}
\caption{Drift $\mu$ and reflection vectors $R^1$ and $R^2$}
\label{fig:drift_reflection}
\end{figure}

The object of study here is the reflected Brownian motion with drift in the quarter plane $\mathbb R_+^2$
\begin{equation}
\label{eq:RBMQP}
     Z(t)=Z_0 + W(t) + \mu t + R L(t),\qquad \forall t\geq 0,
\end{equation}
associated to the triplet $(\Sigma, \mu, R)$, composed of a non-singular covariance matrix, a drift and a reflection matrix, see Figure~\ref{fig:drift_reflection}:
\begin{equation*}
     \Sigma = \left(  \begin{array}{ll} \sigma_{11} & \sigma_{21} \\ \sigma_{12} & \sigma_{22} \end{array} \right),\qquad \mu= \left(  \begin{array}{l} \mu_1 \\  \mu_2  \end{array} \right), \qquad  R= (R^1,R^2)=\left(  \begin{array}{cc} r_{11} & r_{21} \\ r_{12} & r_{22} \end{array} \right).
\end{equation*}     
In Equation~\eqref{eq:RBMQP}, $Z_0$ is any initial point in $\mathbb R_+^2$, the process $(W(t))_{t \geq 0} $ is an unconstrained planar Brownian motion starting from the origin, and for $i\in\{1,2\}$, $L^i(t)$ is a continuous non-decreasing process, that increases only at time $t$ such that $Z^i(t)=0$, namely $\int_{0}^t \mathds{1}_{\{Z^i(s) \ne 0 \}} \mathrm{d} L^i(s)=0$. The columns $R_1$ and $R_2$ represent the directions in which the Brownian motion is pushed when the axes are reached.

The reflected Brownian motion $(Z(t))_{t\geq0}$ associated with $(\Sigma, \mu, R)$ is well defined
, see for instance \cite{Wi95}.
Its stationary distribution exists and is unique if and only if the following (geometric flavored) conditions are satisfied (see, e.g., \cite{harrison_brownian_1987,HoRo93})
\begin{equation}
\label{eq:CNS_ergodic}
     r_{11} > 0, \ \  r_{22} > 0, \ \ r_{11} r_{22} - r_{12} r_{21} > 0,\ \ r_{22} \mu_1 - r_{12}  \mu_2 < 0, \ \ r_{11} \mu_2 - r_{21}  \mu_1 < 0. 
\end{equation}

More that the Brownian motion in the quadrant, all results presented below concern the Brownian motion in two-dimensional cones (by a simple linear transformation of the cones). This is a major difference and interest of the continuous case, which also illustrates that the analytic approach is very well suited to that context.  

\subsection{Functional equation}
\label{subsec:functional_equation_BM}

\subsubsection*{Laplace transform of the stationary distribution}

The continuous analogues of the generating functions are the Laplace transforms. As their discrete counterparts, they characterize the stationary distribution.
Under assumption~\eqref{eq:CNS_ergodic}, that we shall do throughout the manuscript, the stationary distribution is absolutely continuous w.r.t.\ the Lebesgue measure, see \cite{harrison_brownian_1987,Da90}. We denote its density by $\pi(x)=\pi(x_1,x_2)$. Let the Laplace transform of $\pi$ be defined by
\begin{equation*}
     \varphi (\theta) = \mathbb{E}_{\pi} [e^{ \langle \theta \vert  Z\rangle}] = \iint_{{\mathbb R}_+^2} e^{\langle \theta \vert x \rangle} \pi(x) \mathrm{d} x.
\end{equation*}
We further define two finite boundary measures $\nu_1$ and $\nu_2$ with support on the axes, by mean of the formula
\begin{equation*}
     \nu_i (B) = \mathbb{E}_{\pi} \bigg[ \int_0^1 \mathds{1}_{\{Z(t) \in B\}} \mathrm{d}L^i (t)\bigg].
\end{equation*}     
The measures $\nu_i$ are continuous w.r.t.\ the Lebesgue measure by \cite{harrison_brownian_1987}, and may be viewed as boundary invariant measures. We define their moment Laplace transform by
\begin{equation*}
\varphi_2 (\theta_1) 
=\int_{{\mathbb R}_+} e^{\theta_1 x_1} \nu_2(x_1) \mathrm{d} x_1,
\qquad
\varphi_1 (\theta_2) 
=\int_{{\mathbb R}_+} e^{\theta_2 x_2} \nu_1(x_2) \mathrm{d} x_2.
\end{equation*}

\subsubsection*{Functional equation}

There is a functional equation between the Laplace transforms $\phi$, $\phi_1$ and $\phi_2$, see~\eqref{eq:functional_equation}, which is reminiscent of the discrete functional equation \eqref{eq:functional_equation_RW}.

\begin{lem}
The following key functional equation between the Laplace transforms holds
\begin{equation}  
\label{eq:functional_equation}
     -\gamma (\theta) \varphi (\theta) =\gamma_1 (\theta) \varphi_1 (\theta_2) + \gamma_2 (\theta) \varphi_2 (\theta_1),
\end{equation}
where
\begin{align*}
  \begin{cases}
     \phantom{{}_1}\gamma (\theta)= \frac{1}{2} \langle \theta \vert  \sigma \theta \rangle+ \langle \theta \vert  \mu \rangle =\frac{1}{2}(\sigma_{11} \theta_1^2 + \sigma_{22} \theta_2^2 +2\sigma_{12}\theta_1\theta_2)
+
\mu_1\theta_1+\mu_2\theta_2,  \\
     \gamma_1 (\theta)= \langle R^1 \vert  \theta \rangle=r_{11} \theta_1 + r_{21} \theta_2,  \\
     \gamma_2 (\theta)=\langle R^2 \vert  \theta \rangle=r_{12} \theta_1 + r_{22} \theta_2.
  \end{cases}
\end{align*}
\end{lem}
By definition of the Laplace transforms, this equation holds at least for any $\theta=(\theta_1, \theta_2)$ with ${\Re}\,\theta_1\leq 0$ and ${\Re}\, \theta_2\leq 0$. 
The polynomial $\gamma$ in \eqref{eq:functional_equation} is the kernel and is the continuous analogue of the kernel \eqref{eq:kernel_RW} in the discrete case. Polynomials $\gamma_1$ and $\gamma_2$ are the counterparts of $k$ and $\widetilde k$.

\begin{proof}
To show \eqref{eq:functional_equation}, the main idea is to use an identity called a basic adjoint relationship (first proved in \cite{harrison_brownian_1987} in some particular cases, then extended in \cite{DaHa92}), which characterizes the stationary distribution. (It is the continuous analogue of the well-known equation $\pi Q=0$, where $\pi$ is the stationary distribution of a recurrent continuous-time Markov chain with infinitesimal generator $Q$.) This basic adjoint relationship connects the stationary distribution $\pi$ and the corresponding boundary measures $\nu_1$ and $\nu_2$.
We refer to \cite{Fo84,DaMi11} for the details.
\end{proof}

\subsubsection*{Elementary properties of the kernel}

The kernel $\gamma$ in \eqref{eq:functional_equation} can be alternatively written as
\begin{equation}
\label{eq:alternative_expression_kernel_BM}
     \gamma(\theta_1, \theta_2)
     = \widetilde{a}(\theta_2)\theta_1^2+\widetilde{b}(\theta_2)\theta_1+\widetilde{c}(\theta_2)
     =a(\theta_1)\theta_2^2+b(\theta_1)\theta_2+c(\theta_1).
\end{equation}  
The equation $\gamma(\theta_1, \theta_2)= 0$ defines a two-valued algebraic function $\Theta_1(\theta_2)$ by $\gamma(\Theta_1(\theta_2), \theta_2)= 0$, and similarly $\Theta_2(\theta_1)$ such that $\gamma(\theta_1, \Theta_2(\theta_1))= 0$. Expressions of their branches are given by
\begin{equation*}
     \Theta_2^\pm(\theta_1)=\frac{-b(\theta_1)\pm\sqrt{d(\theta_1)}}{2a(\theta_1)},
\end{equation*}
where 
$d (\theta_1) = b^2 (\theta_1) -4a(\theta_1)c(\theta_1)$ is the discriminant.
The polynomial 
$d$ has two zeros, real and of opposite signs; they are denoted by 
$\theta_1^\pm$ and are branch points of the algebraic function $\Theta_2$. In the same way we define $\Theta_1^\pm$ and its branch points $\theta_2^\pm$.

Finally, notice that $d$ is negative on ${\mathbb R} \setminus [\theta_1^-, \theta_1^+]$. Accordingly, the branches $\Theta_2^\pm$ take complex conjugate values on this set. 

\subsection{Statement and resolution of the BVP}
\label{subsec:BVP_BM}

\subsubsection*{An important hyperbola}

\begin{figure}[hbtp]
\centering
\hspace{3mm}\includegraphics[scale=0.45]{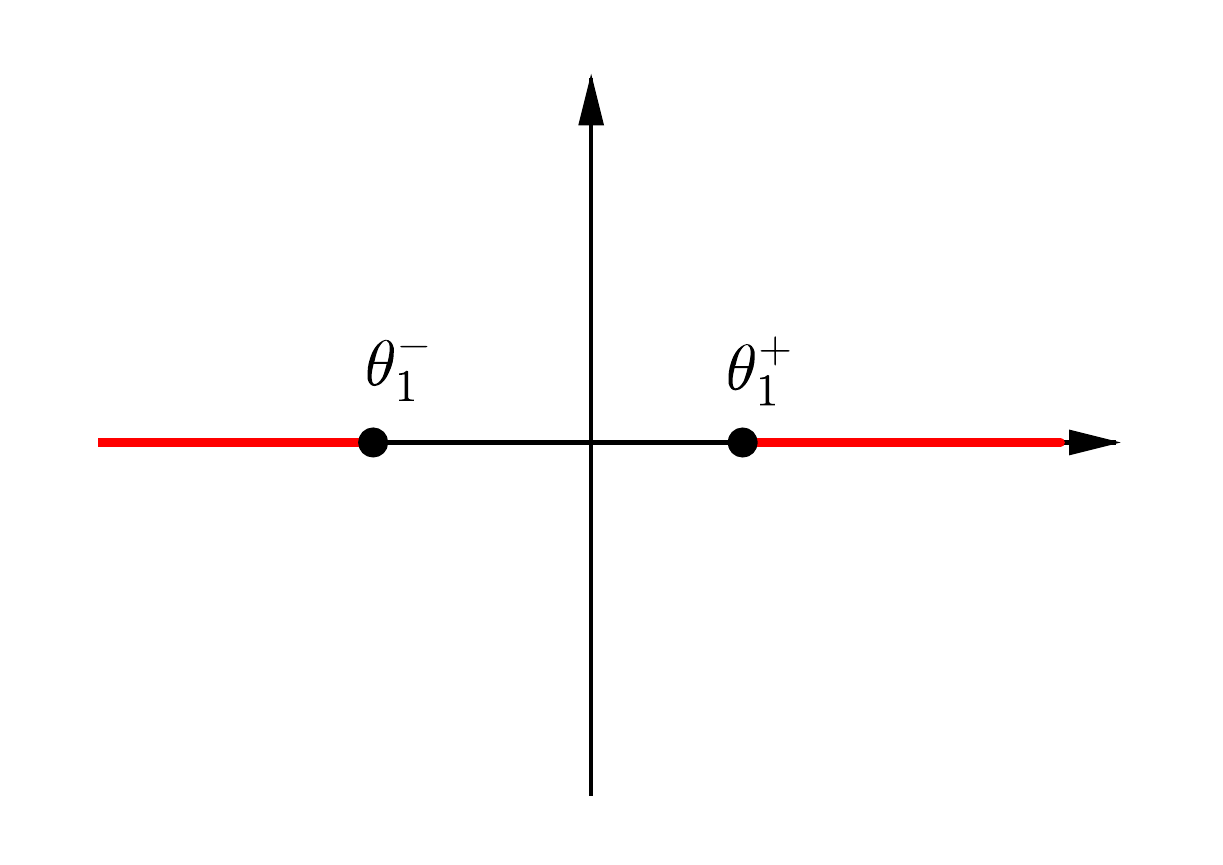}\hspace{8mm}
\includegraphics[scale=0.7]{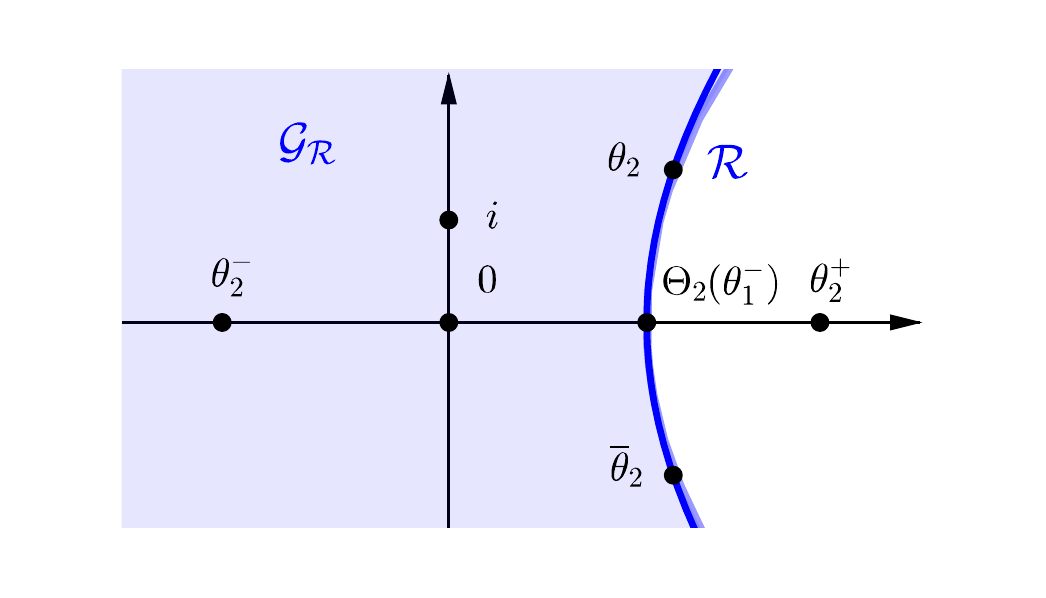}
\caption{Left: the discriminant $d$ has two roots $\theta_1^-$ and $\theta_1^+$ of opposite signs. Right: the curve $\R$ in~\eqref{eq:curve_definition} is symmetric w.r.t.\ the horizontal axis and $\G$ is the domain in blue}
\label{BVPtheta}
\end{figure}

For further use, we need to introduce the curve
\begin{equation}
\label{eq:curve_definition}
     \R =\{\theta_2\in\mathbb C: \gamma(\theta_1,\theta_2)=0 \text{ and } \theta_1\in(-\infty,\theta_1^-)\}=\Theta_2^\pm ((-\infty,\theta_1^-)).
\end{equation}
It is the analogue of the curve $\mathscr{M}$ in Section \ref{subsec:functional_equation_RW}. The curve $\R$ is symmetrical w.r.t.\ the real axis, see Figure~\ref{BVPtheta} (this is a consequence of $d$ being negative on $(-\infty,\theta_1^-)$, see above). Furthermore, it is a {\rm(}branch of a{\rm)} hyperbola by \cite{BaFa87}.
We shall denote by $\G$ the open domain of $\mathbb{C}$ bounded by $\R$ and containing $0$, see Figure~\ref{BVPtheta}. Obviously $\overline{\G}$, the closure of $\G$, is equal to $\G\cup\R$. 

\subsubsection*{BVP for orthogonal reflections}
\label{sec:boundary_cond}

In the case of orthogonal reflections (see Figure \ref{fig:drift_reflection}), $R$ is the identity matrix in~\eqref{eq:RBMQP}, and we have $\gamma_1(\theta_1,\theta_2)=\theta_1$ and $\gamma_2 (\theta_1,\theta_2)=\theta_2$. We set
\begin{equation}
\label{eq:definition_psi}
     \psi_1 (\theta_2)=\frac{1}{\theta_2}\varphi_1 (\theta_2),\qquad \psi_2 (\theta_1)=\frac{1}{\theta_1}\varphi_2 (\theta_1).
\end{equation}

\begin{lem}
\label{prop:BVP}
The function $\psi_1$ in \eqref{eq:definition_psi} satisfies the following BVP:
\begin{enumerate}[label={\rm (\roman{*})},ref={\rm (\roman{*})}]
     \item\label{item:1_BVP} $\psi_1$ is meromorphic on $\G$ with a single pole at $0$, of order $1$ and residue $\varphi_1(0)$, and vanishes at infinity;
     \item\label{item:2_BVP} $\psi_1$ is continuous on $\overline{\G}\setminus \{0\}$ and
\begin{equation}
\label{eq:boundary_condition}
     \psi_1(\overline{\theta_2})=\psi_1({\theta_2}), \qquad \forall \theta_2\in \R.
\end{equation}
\end{enumerate}
\end{lem}

\begin{proof}
The regularity condition of point \ref{item:1_BVP} follows from Theorem \ref{thm:continuation_BM}, which provides a (maximal) meromorphic
continuation of the function.
Let us now consider~\ref{item:2_BVP}. Evaluating the (continued) functional equation \eqref{eq:functional_equation} at $(\theta_1,\Theta_2^\pm(\theta_1))$, we obtain
$
\psi_1(\Theta_2^\pm(\theta_1))+\psi_2(\theta_1)=0,
$
which immediately implies that
\begin{equation}
\label{eq:before_boundary_condition}
     \psi_1(\Theta_2^+(\theta_1))=\psi_1(\Theta_2^-(\theta_1)).
\end{equation}
Choosing $\theta_1\in (-\infty,\theta_1^-)$, the two quantities $\Theta_2^+(\theta_1)$ and $\Theta_2^-(\theta_1)$ are complex conjugate the one of the other, see Section~\ref{subsec:functional_equation_BM}. Equation~\eqref{eq:before_boundary_condition} can then be reformulated as~\eqref{eq:boundary_condition}, using the definition~\eqref{eq:curve_definition} of the curve $\R$.
\end{proof}


The BVP stated in Lemma~\ref{prop:BVP} is called a homogeneous BVP with shift (the shift stands here for the complex conjugation, but the theory applies to more general shifts, see \cite{Li00}). It has a simpler form than the BVP in Lemma~\ref{lem:BVP_RW} for the discrete case, because there is no inhomogeneous term (as $\pi_{0,0}$) and also because in the coefficients in front of the unknowns there is no algebraic function (as $Y_0$) involved. Due to its particularly simple form, we can solve it in an explicit way, using the two following steps: 
\begin{itemize}
     \item Using a certain conformal mapping $w$ (to be introduced below), we can construct a particular solution to the BVP of Lemma~\ref{prop:BVP}.
     \item The solution to the BVP of Lemma~\ref{prop:BVP} is unique (see the invariant Lemma 2 in \cite[Section 10.2]{Li00}). In other words, two different solutions must coincide, and the explicit solution constructed above must be the function $\psi_1$.
\end{itemize}
In \cite{FrRa16} it is explained that the above method may be viewed as a variation of Tutte's invariant approach, first introduced by Tutte for solving a functional equation arising in the enumeration of properly colored triangulations, see \cite{Tutte-95}.

The function $w$ glues together the upper and lower parts of the hyperbola $\R$.
There are at least two ways to find such a $w$. First, it turns out that in the literature there exist expressions for conformal gluing functions for relatively simple curves as hyperbolas, see \cite[Equation (4.6)]{BaFa87}. Here (based on \cite{FrRa16}), we use instead the Riemann sphere $\s$, as we will see in Section \ref{subsec:Riemann_surface_BM}. Indeed, many technical aspects (and in particular finding the conformal mapping) happen to be quite simpler on that surface. 

We will deduce from Section \ref{subsec:Riemann_surface_BM} that function $w$ can be expressed in terms of the generalized Chebyshev polynomial 
\begin{equation*}
T_a(x)  =\cos (a\arccos (x))=\frac{1}{2} \Big\{\big(x+\sqrt{x^2-1}\big)^a+\big(x-\sqrt{x^2-1}\big)^a\Big\}
\end{equation*} as follows:
\begin{equation}
\label{eq:expression_CGF_BM}
     {w} (\theta_2)=T_{\frac{\pi}{\beta}}\bigg(-\frac{2\theta_2-(\theta_2^++\theta_2^-)}{\theta_2^+-\theta_2^-}\bigg),
\end{equation}
where we have noted
\begin{equation}
\label{eq:definition_beta}
     \beta=\arccos {-\frac{\sigma_{12}}{\sqrt{\sigma_{11}\sigma_{22}}}}.
\end{equation}
In the case of orthogonal reflection, this methods leads to the main result of \cite{FrRa16}, which is:
\begin{thm} 
\label{thm:main}
Let $R$ be the identity matrix in~\eqref{eq:RBMQP}. The Laplace transform $\phi_1$ is equal to 
\begin{equation*}
     \phi_1 (\theta_2)= \frac{-\mu_1 {w}'(0)}{{w}(\theta_2)-{w}(0)}\theta_2.
\end{equation*}
\end{thm}

\subsubsection*{Statement of the BVP in the general case}
\label{subsec:general_case}

We would like to close Section~\ref{subsec:BVP_BM} by stating the BVP in the case of arbitrary reflections (non-necessarily orthogonal). 
Let us define for $\theta_2\in \R$
\begin{equation*}
     G(\theta_2)=\frac{\gamma_1}{\gamma_2}(\Theta_1^-(\theta_2),\theta_2)\frac{\gamma_2}{\gamma_1}(\Theta_1^-(\theta_2),\overline{\theta_2}).
\end{equation*}     
Similarly to  Lemma~\ref{prop:BVP}, there is the following result:
\begin{lem}
\label{prop:BVP_general}
The function $\phi_1$ satisfies the following BVP:
\begin{enumerate}[label={\rm (\roman{*})},ref={\rm (\roman{*})}]
     \item\label{item:1_BVP_general} $\phi_1$ is meromorphic on $\G$ with at most one pole $p$ of order $1$ and is bounded at infinity;   
     \item\label{item:2_BVP_general} $\phi_1$ is continuous on $\overline{\G}\setminus \{p\}$ and
\begin{equation}
\label{eq:boundary_condition_general}
     \phi_1(\overline{\theta_2})=G(\theta_2)\phi_1({\theta_2}), \qquad \forall \theta_2\in \R.
\end{equation}
\end{enumerate}
\end{lem} 
Due to the presence of the function $G\neq1$ in \eqref{eq:boundary_condition_general}, this BVP (still homogeneous with shift) is more complicated than the one encountered in Lemma~\ref{prop:BVP} and cannot be solved thanks to an invariant lemma. Instead, the resolution is less combinatorial and far more technical, and the solution should be expressed in terms of both Cauchy integrals and the conformal mapping $w$ of Theorem~\ref{thm:main}. This will be achieved in a future work.

\subsection{Asymptotics of the stationary probabilities}
\label{subsec:asymp_BM}

\subsubsection*{Overview}

Let  $\Pi$ be a random vector  that has the stationary distribution of the reflected Brownian motion.
\cite{DaMi11} obtain the following asymptotic result:
for a given directional vector $c \in {\mathbb{R}}_+^2$ they find (up to a multiplicative constant) a function $f_c(x)$ such that
\begin{equation*}
     \lim_{ x \to \infty} \frac{\mathbb{P}[ \langle c  \vert \Pi \rangle \geq x]}{f_c(x)}=1.
\end{equation*}
In \cite{FrKu16} we solve a harder problem arisen in \cite[\S 8]{DaMi11}, namely computing the asymptotics of $\mathbb{P}[\Pi \in  x c+B]$ as $x\to \infty$, 
where $ c \in {\mathbb{R}}_+^2$ is any directional vector and $B \subset {\mathbb{R}}_+^2$ any compact subset. Furthermore, we are able to find the full asymptotic expansion  of the density $\pi(x_1,x_2)$ of $\Pi$  as $x_1,x_2 \to \infty$  and $x_2/x_1 \to  {\rm \tan } (\alpha)$, for any given angle  $\alpha \in (0,\pi/2)$.

\subsubsection*{Main results}

\begin{figure}[hbtp]
\centering 
\includegraphics[scale=0.7]{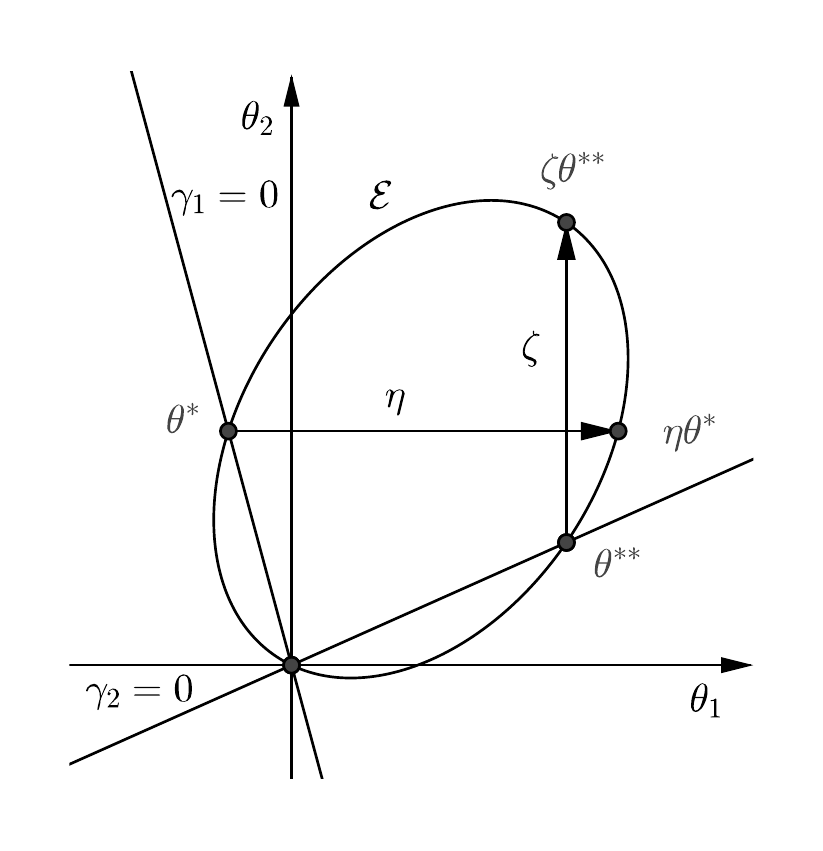} 
\includegraphics[scale=0.7]{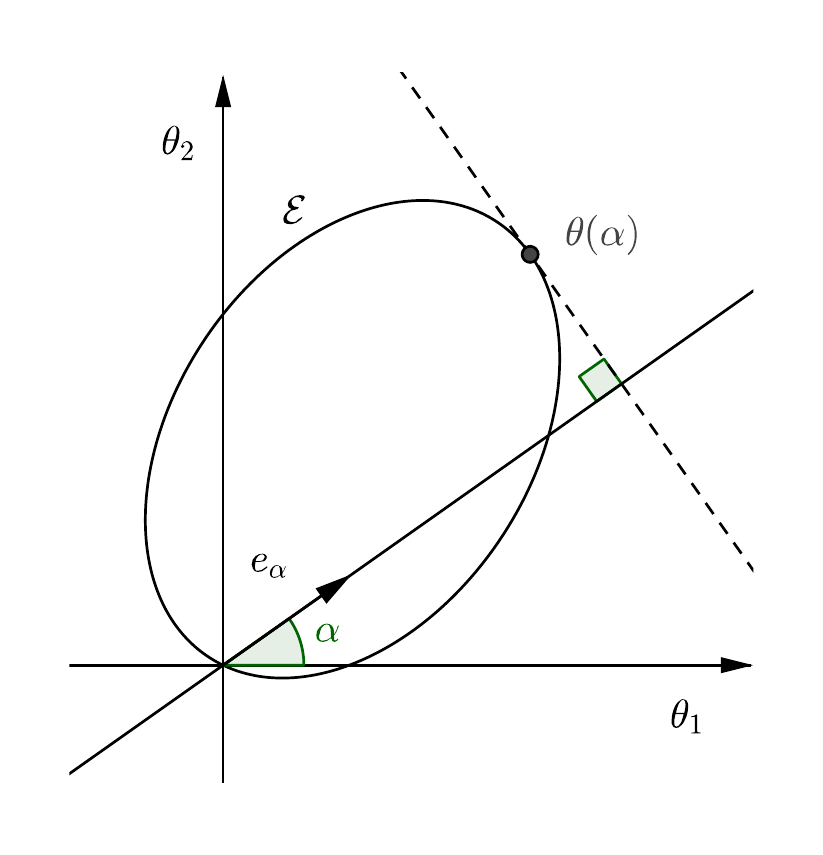} 
\caption{Left: representation of the ellipse $\mathscr{E}$, straight lines $\{\gamma_1(\theta)=0\}$, $\{\gamma_2(\theta)= 0\}$, and points 
$\theta^{*}$, $\theta^{**}$, $\eta \theta^*$ and $\zeta \theta ^{**}$. Right: geometric interpretation of the point $\theta(\alpha)$ in \eqref{thetaalpha} on $\mathscr{E}$}
 \label{ellipse}
 \label{pointpolgeom}
\end{figure}

First we need to introduce some notations. The equation $\gamma(\theta)=0$ determines an ellipse $\mathscr{E}$ on ${\mathbb R}^2$ passing through the origin,
see Figure \ref{ellipse}.
Here we restrict ourselves to the case
$\mu_1<0$ and $\mu_2<0$,
although our methods can be applied without additional difficulty to other cases.
For a given angle $\alpha \in [0,\pi/2]$, let us define the point $\theta(\alpha)$ on the ellipse ${\mathscr E}$ by
\begin{equation}
\label{thetaalpha}
\theta(\alpha)={\rm argmax}_{\theta \in {\mathcal E}}\langle \theta \vert e_\alpha  \rangle, \qquad \hbox{where }e_\alpha=(\cos \alpha, \sin \alpha).
\end{equation}
The coordinates of $\theta(\alpha)$ can be given explicitly.
One can also construct $\theta(\alpha)$ geometrically
as on Figure~\ref{pointpolgeom}.

Secondly, consider the straight lines $\{\gamma_1(\theta)=0\}$ and $\{\gamma_2(\theta)=0\}$, depending on the reflection matrix $R$ only.
They cross  the ellipse ${\mathscr E}$ at the origin. The  line $\{\gamma_1(\theta)=0\}$ (resp.\ $\{\gamma_2(\theta)=0\}$)  intersects the ellipse at a second point called $\theta^{*}$ (resp.\ $\theta^{**}$).
To present our results, we need to define the images on ${\mathscr E}$ of these points via the so-called Galois automorphisms $\zeta$ and $ \eta$, to be introduced in Section \ref{subsec:Riemann_surface_BM}.
Namely, for the point $\theta^*=(\theta_1^*,\theta_2^*) \in {\mathscr E}$ there exists
a unique point $\eta \theta^*=(\eta \theta_1^*, \theta_2^*) \in {\mathscr E}$ with the same second coordinate.
Likewise, there exists a unique point  $\zeta \theta^{**}=( \theta_1^{**},  \zeta \theta_2^{**}) \in {\mathscr E}$  with the same first coordinate as $\theta^{**}=(\theta_1^{**},\theta_2^{**}) \in {\mathscr E}$.
Points  $\theta^*$, $\theta^{**}$,  $\eta \theta^*$ and $\zeta \theta^{**}$ are pictured  on Figure \ref{ellipse}.
Their coordinates can be made explicit.

Similarly to the discrete case, we introduce the set of parameters
\begin{equation*}
 \mathscr{Q}_{--}=\big\{((\Sigma,\mu,R), \alpha): \gamma_1(\eta\theta(\alpha)) <0 \text{ and } \gamma_2(\zeta\theta(\alpha)) <0 \big\}\nonumber
\end{equation*}
and $\mathscr{Q}_{+-}$, $\mathscr{Q}_{-+}$ and $\mathscr{Q}_{++}$ accordingly.
The following theorem provides the main term in the asymptotic expansion of $\pi(r\cos \alpha, r \sin \alpha)$.
\begin{thm}
   Let $(x,y)=(r \cos \alpha,r \sin \alpha)$ with $\alpha \in (0, \pi/2)$. We assume that $\theta (\alpha) \in \mathbb{R}_+^2$.
Then as $r \to \infty$  we have
\begin{equation}
\pi (r e_\alpha) = 
  (1+o(1))\cdot  \left\{
\begin{array}{ll}
\frac{C_0}{\sqrt{r}}e^{-r \langle e_\alpha \vert\theta(\alpha) \rangle}  & \hbox{ in } \mathscr{Q}_{--},\\
C_1e^{-r \langle e_\alpha \vert\eta\theta^{*}\rangle} & \hbox{ in } \mathscr{Q}_{+-},\\
C_2 e^{-r \langle e_\alpha \vert\zeta\theta^{**}\rangle}& \hbox{ in } \mathscr{Q}_{-+},\\
C_1e^{-r \langle e_\alpha \vert\eta\theta^{*}\rangle}  +
 C_2e^{-r \langle e_\alpha \vert\zeta\theta^{**}\rangle} & \hbox{ in } \mathscr{Q}_{++},
\end{array}
\right.
\end{equation}
   where $C_0$, $C_1$ and $C_2$ are constants that can be expressed in terms of functions $\phi_1$ and $\phi_2$ and the parameters.
   \label{thmmain}
\end{thm}
In \cite{FrKu16} the constants mentioned in Theorem \ref{thmmain} are specified in terms of functions $\varphi_1$ and $\varphi_2$. But these functions
are for now unknown. 
As we explained in Section \ref{subsec:BVP_BM}, in a next work we are going to obtain $\phi_1$ and $\phi_2$ as solutions of BVP, thereby determining the constants in Theorem~\ref{thmmain}.

\begin{proof}[of the key step of Theorem \ref{thmmain}]
 Theorem \ref{thmmain} is proven in \cite{FrKu16}.
The first step consists in continuing meromorphically the functions $\phi_1$ and $\phi_2$ on $\mathbb{C}\setminus [\theta_2^+,\infty)$ or on the Riemann surface $\s$, see Section \ref{subsec:Riemann_surface_BM}.
Then by the  functional equation \eqref{eq:functional_equation}  and the  inversion formula of  Laplace transform
 (we refer to \cite{doetsch_introduction_1974} and
\cite{brychkov_multidimensional_1992}), 
       the density $\pi(x_1,x_2)$ can be represented as a double integral. Using standard computations from complex analysis, we are able to reduce it to a sum of single integrals. We obtain the following (with the notation \eqref{eq:alternative_expression_kernel_BM}):
\begin{align*}
\pi(x_1,x_2)
&=\frac{-1}{(2\pi i)^2} \int_{-i\infty}^{i\infty} \int_{-i\infty}^{i\infty}
e^{-x_1\theta_1-x_2\theta_2}
\frac{\gamma_1(\theta) \varphi_1(\theta_2)+\gamma_2(\theta)\varphi_2(\theta_1)}{\gamma(\theta)}
\mathrm{d} \theta_1 \mathrm{d} \theta_2 
\\
&= \frac{1}{2\pi i} \int_{-i\infty}^{i\infty}
\varphi_2(\theta_1) \gamma_2(\theta_1, \Theta_2^+(\theta_1)) e^{-x_1\theta_1-x_2\Theta_2^+(\theta_1)}
\frac{\mathrm{d} \theta_1}{\sqrt{  d(\theta_1) }} 
\\ &+\frac{1}{2\pi i} \int_{-i\infty}^{i\infty}
\varphi_1(\theta_2) \gamma_1(\Theta_1^+(\theta_2), \theta_2) e^{-x_1\Theta_1^+(\theta_2)-x_2\theta_2}
\frac{\mathrm{d} \theta_2}{\sqrt{  \widetilde d (\theta_2) }}.
\end{align*}
These integrals are typical to apply the saddle point method, see \citet{fedoryuk_asymptotic_1989}. 
The coordinates of the saddle point are the critical points of the functions
\begin{equation*}
\cos(\alpha)\theta_1+\sin(\alpha)\Theta_2^+(\theta_1)\quad \text{ and }\quad\cos(\alpha)\Theta_1^+(\theta_2)+\sin(\alpha)\theta_2. 
\end{equation*}
It is the point $\theta(\alpha)$. Then we have to shift the integration contour up to new contours which coincide with the steepest-descent contour near the saddle point.
When we shift the contours we have to take into account the poles of the integrands and their residues. The asymptotics will be determined by the pole if we cross a pole when we shift the contour and by the saddle point otherwise.
\end{proof}

\subsection{Riemann surface and related facts}
\label{subsec:Riemann_surface_BM}

\subsubsection*{Riemann surface}
\label{subsubsec:riemannsurface}

The Riemann surface 
\begin{equation*}
     \s=\{(\theta_1,\theta_2)\in\mathbb{C}^2 : \gamma(\theta_1,\theta_2)=0\}
\end{equation*}
may be viewed as the set of zeros of the kernel (equivalently, it is the Riemann surface of the algebraic functions $\Theta_2$ and $\Theta_1$). Due to the degree of $\gamma$, the surface $\s$ has genus $0$ and is a Riemann sphere, i.e., homeomorphic to $\mathbb{C}\cup\{\infty\}$, see \cite{FrKu16}.
It admits a very useful rational parametrization, 
given by
\begin{equation}
\label{eq:formulas_uniformization}
     \theta_1(s)=\displaystyle\frac{\theta_1^{+}+\theta_1^{-}}{2}+\frac{\theta_1^{+}-\theta_1^{-}}{4} \left(s+\frac{1}{s}\right),\qquad
      \theta_2(s)=\displaystyle\frac{\theta_2^{+}+\theta_2^{-}}{2}+\frac{\theta_2^{+}-\theta_2^{-}}{4} \left(\frac{s}{e^{i\beta}}+\frac{e^{i\beta}}{s}\right),
\end{equation}
with $\beta$ as in \eqref{eq:definition_beta}. The equation $\gamma(\theta_1(s),\theta_2(s))=0$ holds and $\s= \{(\theta_1(s),\theta_2(s)): s\in\mathbb{C}\cup\{\infty\} \}$. 

\subsubsection*{Group of the process}

We finally introduce the notion of group of the model, similar to the notion of group of the walk in the discrete setting (see \cite{Ma72,FaIaMa99,BMMi10}).
This group $\langle\zeta,\eta\rangle$ is generated by $\zeta$ and $\eta$, given by (with the notation \eqref{eq:alternative_expression_kernel_BM})
\begin{equation*}
     \zeta(\theta_1,\theta_2)=\left(\theta_1,\frac{c(\theta_1)}{a(\theta_1)}\frac{1}{\theta_2}\right),\qquad
     \eta(\theta_1,\theta_2)=\left(\frac{\widetilde c(\theta_2)}{\widetilde a(\theta_2)}\frac{1}{\theta_1},\theta_2\right).
\end{equation*}
By construction, the generators satisfy $\gamma(\zeta(\theta_1,\theta_2))=\gamma(\eta(\theta_1,\theta_2))=0$ as soon as $\gamma(\theta_1,\theta_2)=0$. In other words, there are (covering) automorphisms of the surface $\s$. Since $\zeta^2=\eta^2=1$, the group $\langle\zeta,\eta\rangle$ is a dihedral group, which is finite if and only if the element $\zeta\circ\eta$ (or equivalently $\eta\circ\zeta$) has finite order.

\subsubsection*{Algebraic nature of the Laplace transforms}
With the above definition, it is not clear how to see if the group is finite, nor to see it its finiteness would have any implication on the problem. In fact, we have, with $\beta$ defined in \eqref{eq:definition_beta}:
\begin{lem}
\label{prop:group_finite}
The group $\langle\zeta,\eta\rangle$ is finite if and only if $\pi/\beta\in\mathbb Q$.
\end{lem}
The proof of Lemma~\ref{prop:group_finite} is simple, once the elements $\zeta$ and $\eta$ have been lifted and reformulated on the sphere $\s$:
\begin{equation*}
     \zeta(s)=\frac{1}{s},\qquad \eta(s)=\frac{e^{2i\beta}}{s}.
\end{equation*}
These transformations leave invariant $\theta_1(s)$ and $\theta_2(s)$, respectively, see \eqref{eq:formulas_uniformization}. In particular, we have the following result (see \cite{FrRa16}), which connects the nature of the solution of the BVP to the finiteness of the group. Such a result holds for discrete walks, see our Theorem \ref{thm:nature_solution_RW} and \cite{BMMi10,BeBMRa15}. 
\begin{thm}
\label{prop:nature}
The solution $\phi_1$ given in Theorem~\ref{thm:main} and the conformal gluing function $w$ in \eqref{eq:expression_CGF_BM} are algebraic if and only if the group $\langle\zeta,\eta\rangle$ is finite.
\end{thm}

\subsubsection*{Conformal mapping}
The conformal gluing function $w$ introduced in Section \ref{subsec:BVP_BM} may be lifted on $\s$. In fact its expression is even simpler using the parametrization of $\s$. We show in \cite{FrRa16} that
\begin{equation}
\label{eq:expression_W}
     w(\theta_2(s))=-\frac{i}{2}\big\{(-s)^{\frac{\pi}{\beta}}+(-s)^{-\frac{\pi}{\beta}}\big\}=-\frac{i}{2}\big\{e^{\frac{\pi}{\beta}\log (-s)}+e^{-\frac{\pi}{\beta}\log (-s)}\big\},
\end{equation}
where we make use of the principal determination of the logarithm.

\subsubsection*{Continuation of the Laplace transforms}
\label{subsubsec:continuation}

To establish the BVP, we have stated a boundary condition for the functions $\phi_1$ and $\phi_2$, on curves which lie outside their natural domains of definition (the half-plane with negative real-part), see Figure \ref{BVPtheta}. In the same way, in the asymptotic study we use the steepest descent method on some curves outside of the initial domain of definition. We therefore need to extend the domain of definition of the Laplace transforms.
\begin{thm}
\label{thm:continuation_BM}
The function $\phi_1$ can be continued meromorphically on the cut plane $\mathbb{C}\setminus [\theta_2^+,\infty)$.
\end{thm}

\begin{proof}
The first step is to continue meromorphically $\phi_1(\theta_2)$ to the open and simply connected set
$
     \{\theta_2\in\mathbb{C} : \Re
    \,
     \theta_2\leq
     0 \text{ or } \Re\, \Theta_1^-(\theta_2) <0\},
$
by setting
\begin{equation*}
     \phi_1(\theta_2)=\frac{\gamma_2}{\gamma_1}(\Theta_1^-(\theta_2),\theta_2)\phi_2(\Theta_1^-(\theta_2)).
\end{equation*}
This is immediate (see \cite{FrKu16} for the details). 
It is then possible to pursue the extension to the whole $\s$ using the invariance properties by the automorphisms $\zeta$ and $\eta$ satisfied by the lifted Laplace transforms on $\s$.
\end{proof}

\acknowledgements
We acknowledge support from the projet \href{http://www.fdpoisson.fr/madaca/}{MADACA} of the R\'egion Centre and from Simon Fraser University (British Columbia, Canada). We thank our three AofA referees for useful comments.


\footnotesize

\bibliographystyle{apalike} 
\bibliography{biblio} 
\normalsize
\end{document}